\documentclass{amsart}

\usepackage{amsmath, amsthm}
\usepackage{amssymb, latexsym}
\usepackage{amsfonts}

\usepackage[dvipsnames]{xcolor}

    \usepackage{lineno}
\newtheorem{Theorem}[equation]{Theorem}
\newtheorem{theorem}[equation]{Theorem}

\newtheorem{Example}[equation]{Example}

\newtheorem{Remark}[equation]{Remark}
\newtheorem{Remarks}[equation]{Remarks}

\theoremstyle{definition}

\newcommand{\exref}[1]{Ex\-am\-ple \ref{#1}}
\newcommand{\thmref}[1]{\rm Theorem~\ref{#1}}

\newcommand{\eqnref}[1]{\rm(\ref{#1})}

\newcommand{\beql}[1]{\begin{equation}\label{#1}}
\newcommand{\eeq} {\end{equation}}

    \font\Aaa=msam10

\font\Bbb=msbm10

\newcommand\F{\hbox{\Bbb F}}

\font\Aaa=msam10

\def\qed{\hbox{~~\Aaa\char'003}}

\font\Bbb=msbm10

\def\F{\hbox{\Bbb F}}

\numberwithin{equation}{section}

\DeclareMathOperator{\Aut}{Aut}

\newcommand\G{\Gamma }
\newcommand\D{{ \Delta }}
\newcommand\z{{ \zeta }}

        \def\PSL{{\rm PSL}}

        \def\Sp{{\rm Sp}}
        
        \def\PG{{\rm PG}}

        \def\<{{\langle}}
        \def\>{{\rangle}}

        \font\Aaa=msam10

\font\Aaa=msam10

\font\Bbb=msbm10

\def\F{\hbox{\Bbb F}}
\font\Bbbsmall=msbm7

\def\Zsmall{\hbox{\Bbbsmall Z}}

\def\div{ \kern-.5pt\hbox{\big |} }
\def\ndiv{ {\not\kern-.5pt\hbox{\big |}\,} }
\def\ndivv{ {\not\kern+1.5pt\hbox{$\mid$}\,} }

\def\Aut{{\rm Aut}}

\def\dim{{\rm dim}}

\def\col{\colon\!}

\def\B{^2\kern-.8pt B}
\def\G{^2\kern-.8pt G}
\def\EH{^2\kern-.8pt\hat  E}
\def\E{^2\kern-.8pt E}
\def\D{^3\kern-1pt D}
\def\FF{^2\kern-.8pt F}

\newdimen\refcodesize
\newbox\seriesbox
\def\proof{\noindent {\bf Proof.~}}

\def\PSL  {{\rm PSL }}

\def\Sp{{\rm Sp}}

\def\2{^{(2)}}

\def\ffiti{\big((F_i)_0^{n},(\zeta_i)_0^{n-1}\big)}
\def\ffisi'{\big((F_i')_0^{n'},(\zeta'_i)_0^{n'-1}\big)}

\def\fiti{\big((F_i)_0^{n},(\zeta_i)_0^{n}\big)}

\def\fisi'{\big((F_i')_0^{m},(\zeta'_i)_0^{m}\big)}
\def\fisim'{\big((F_i')_0^{n'},(\zeta'_i)_0^{n'}\big)}

\def\Aut{{\rm Aut}}
\def\2{^{(2)}}

\def\c{{\rm C}}

\def\sss{{\mathcal S}}

\def\c{{\rm C}}

\def\vartheta{\gamma}
\def\tc{{\widetilde\c}}

\def\2{^{(2)}}
\def\Aut{{\rm Aut}}

\DeclareRobustCommand{\SkipTocEntry}[4]{}
\usepackage{lineno}
 
\begin{document} 

\title[Elliptic  spreads]
{Cyclic elliptic  spreads
 }
 
       \author{W. M. Kantor}
       \address{U. of Oregon,
       Eugene, OR 97403
        \      and  \
       Northeastern U., Boston, MA 02115}
       \email{kantor@uoregon.edu}

\begin{abstract}
\vspace{-4pt} 

Moderately large numbers of   transitive elliptic spreads   are constructed in
 characteristic $2$ and   dimension $\equiv $ 2 (mod 4).
 \end{abstract}

\maketitle
\vspace{-24pt}
   \section{Introduction}
   
If $V$ is an orthogonal  vector space,~an 
 \emph{orthogonal    spread}
is   a set of maximal totally singular  subspaces  such that   every nonzero   
    singular  vector is in a unique member of  the set.    
    Throughout this note  $q$ will denote a power of 2. 
   Dillon \cite{Dillon}, and later also   Dye \cite{Dye},  showed 
   that $O^-(2m,q)$-spaces and $O(2m-1,q)$-spaces  have orthogonal spreads 
   for all $m\ge2$,
   and that   $O^+(2m,q)$-spaces have orthogonal spreads if and only if  
   $m$ is  even.
   Moreover, the    $O^+(2m,q)$-spreads and $O(2m-1,q)$-spreads 
   constructed in  
    \cite{Dillon,Dye}
    are permuted   3-transitively by isometry  groups isomorphic to 
    $\PSL(2,q^{m-1})$, which   is  not the case for  the   
    $O^-(2m,q)$-spreads in  \cite{Dillon,Dye} when  $m>2$  (see \exref{from desarguesian} below).
   The fact that the various orthogonal spreads  in  
   those papers   admit transitive cyclic groups of isometries is implicit in their constructions (see \cite{Ka,KW1} and \exref{from desarguesian}).

  Much of the small literature on elliptic spreads has focused on existence in large dimensions (citing \cite{Dye}, cf.  \exref{symplectic to orthogonal} below)  or examples in dimension 6 behaving ``nicely'' (as in \cite {CK,CE}).
  This note is concerned with large-dimensional examples,
  using 
  \cite{KW1} to obtain   many $O^-(2m,q)$-spreads
  admitting  transitive cyclic groups of isometries 
  (although ``many'' is rather small compared to the numbers obtained in 
  related research on semifields \cite{KW2}):
  
    \Theorem
  \label{Main Theorem}
  Let $q$ be a power of $\,2$.    Let
$m, m_1,\ldots,m_{n-1} $ be any sequence~of 
$n \ge 1$ distinct divisors $> 1$ of     an odd
  integer   $m$  
   with
each  divisible by the next.  Then there are more than
$\big( \prod_1^{n-1}(q^{m_i}+1)\big)/{2m_1}\log q $ pairwise inequivalent 
$O^-(2m,q)$-spreads each of which is permuted transitively by a cyclic group of isometries.

  \medskip
  \rm
   We will see that this result is closely tied to results  in  \cite{KW1}.
   Such a sequence $(m_i)$ can be obtained from the prime factorization 
   $m=\prod_1^np_j$ of $m$ by setting $m_i=\prod_{i+1}^n p_j$.  In particular, we obtain more than   $q^{m/p}\big/2(m/p)\log q$ orthogonal spreads when
   $ p$ is the  smallest prime dividing $m$.
   
   The desirability of having examples of  transitive orthogonal spreads can be seen from \cite{BP}.
   There are analogous results  in \cite{KW1} producing many $O^+(2m,q)$-spreads and $O(2m-1,q)$-spreads.

  \section{Background}
  
   We refer to  \cite{Taylor} for the standard properties of the orthogonal and symplectic vector spaces  used here.  We name  geometries  using their  isometry groups.  
We will be concerned with singular vectors  and totally singular (t.s.) 
subspaces of  orthogonal spaces, and totally isotropic (t.i.) subspaces of   symplectic spaces.  
In characteristic 2, an orthogonal   vector space    
is also 
a  symplectic space, and   t.s. subspaces are also t.i.  subspaces.

Orthogonal spreads were defined earlier.   For the $O^-(2m,q)$-spaces considered here, an orthogonal spread consists of
    $q^m+1$ t.s. $m-1$-spaces, and is also called an 
    \emph{elliptic spread}.
A      \emph{symplectic     spread}
  is   a set of maximal t.i.   subspaces  of a symplectic space such that   every nonzero vector   is in a unique member of  the set.  

Two orthogonal
or  symplectic  spreads are 
\emph{equivalent} if there is an 
isomorphism of the underlying  
orthogonal
or  symplectic   geometries sending one  spread to the other.  
   The automorphism group of  an orthogonal
or  symplectic    spread    is 
   the group of such isomorphisms of the spread   with itself. 

    \section{Examples}
  \label{Examples}
Let  $F\2=\F_{q^{2m}}\supset F=\F_{q^{m}}\supset K=\F_{q } $
for $m\ge2$, with involutory field  automorphism
$x\mapsto \bar x$.
Let $W:=\ker T $  for the   trace map   $T\col F\to K$.  

The quadratic form  $Q(x):=T(x \bar x)$ turns   $F\2$   into an orthogonal $K$-space with  associated alternating bilinear form    
$(x,y):=
T(x \bar y  +\bar x y)$.   The subspace $W$ is t.s.:  if $w\in W$ then 
$Q(w)=T(w \bar w)=T(w)^2=0$.   

Write $\c:= \{\theta\in F\2\mid \theta \bar\theta=1 \}$, so that
$F\2{}^* = F^*\times \c$.  If $\theta\in\c$ then 
$\widetilde\theta\col x\mapsto x\theta $ defines an isometry
of  the orthogonal space $F\2 $.  The group  $\tc$  of these 
$q^m+1$ isometries
 makes it clear that $F\2 $  is  an
 $O^-(2m,q)$-space: its singular points form an elliptic quadric.
Moreover,  $W$ is a maximal t.s. subspace since  $\dim W=m-1$.%

  \Example
  \label{from desarguesian}
  \rm
   The set $F^\tc$ is the usual desarguesian spread.  Moreover, $W^\tc$ is an elliptic spread  (compare \cite[p.~188]{Dye})
   permuted transitively by 
 the  cyclic    isometry group    $\tc$:  the members of $W^\tc$  intersect pairwise in 0, and its union has
   size
    $1+(q^m +1)(q^{m-1}-1)$, which is the number of singular vectors.
In view of \cite[Theorem~1.1]{BP}, if $m>2$ then the automorphism group of 
$W^\tc$   normalizes  $\tc$.

      In dimension $2m=6$  the Klein correspondence produces 
      an   ovoid  from $W^\tc$ that is constructed in an entirely different manner 
      and  studied in detail in  \cite{CK}.  These examples are equivalent  by 
     \cite[Theorem~1.1]{CK}, and the
  automorphism group appears in 
\cite[Theorem~3.1]{CK}.

  \Example\rm
  \label{symplectic to orthogonal}
   Let
  $\Sigma$ be any symplectic spread in an $\Sp(2m,q)$-space $V$.
  There are many  different ways that $V$ can be equipped with 
  the structure of an 
  $O^-(2m,q)$-space producing the given symplectic structure.    
  Choose one of them.  We emphasize that 
  this orthogonal structure is not related to   $\Sigma$.
  
     If $X\in \Sigma$ let 
  $X'$ be the set of singular vectors in $X$.     Then $X'$ is a t.s. 
  $(m -1)$-space (since the quadratic form $Q|_X\col X\to K$ is semilinear on the t.i. subspace $X$),  and 
   $ \{ X' \mid X\in \Sigma \}$
  is an orthogonal spread of the orthogonal space $V$.
  
  This presumably produces reasonably large numbers of inequivalent elliptic spreads starting from a single symplectic spread.  
  However,  in the next section   we will proceed differently, using inequivalent symplectic spreads in order to arrange  that all types of spreads
 are  permuted transitively by the same cyclic group.
 
 Most of the {\em known}  symplectic spreads in characteristic 2 are the ones in  
 \exref{from desarguesian} or are obtained by a process described  in \cite{Ka,KW1,KW2}; examples are in the next section.
 The only other types  arise from the Suzuki groups \cite{Lu},
   are in $\Sp(4m,2^e)$-spaces for odd   $m $  and  $e$,
 and  again produce elliptic spreads.

    \section{Moderately large numbers of elliptic spreads}
The following notation is based on  \cite{KW1}.
Let $F\2\supset F=F_0\supset\cdots\supset F_{n} =\F_q$ be a tower of fields with
$m:=[F\col F_{n }]$ odd  and corresponding trace maps  
$T_i\col  F\to F_i$.  Write
$W_i:=\ker (T_{i+1}\big|_{F_i})$.
Since $ m$ is odd, $T_{n }(1)=1$ and 
\vspace{-1pt}
$T_{n } T_i(x)= T_{n }(x)$ for all $i$ and all $x\in F$. 
For each $i$  let $F_i\2$ be the subfield of $F\2$ of degree 2 over
$F_i$.

View $V:=F\2$ as
an    $O^-(2m,q)$-space   with  
associated quadratic  form   
$Q_{n}(x):=T_{n }(x\overline{x} ) $ and  
 alternating form   
$ (x,y )_{n} :=T_{n }(x\overline{y}+\overline{x}y)$.   
Let $ \c$ and    $\tc$ be as before.

Let  $\zeta_0=1$ and 
 $1\ne \zeta_i\in\c\cap F_i\2$  for  
$i\geq1$.
Write $\gamma_i:=\prod_0^i\zeta_j, 0\le i\le m-1,$~and   
\begin{equation}
\label{symplectic spread}
\sss\fiti:=\big\{
\big(\sum _0^{n-1}W_i\vartheta_i+F_{n }\vartheta_{n }\big)
\theta \mid \theta\in\c\big\}.
\end{equation}
By \cite[Theorems ~5.1~and~5.2]{KW1}
 (or  \cite[pp.~565-617]{JJB}),
{\em $\sss\fiti$ is a symplectic spread of the
$F_{n }$-space $V,$ and
$\tc$ acts transitively on this spread.}
Moreover, {\em for  $\fiti$ and $\fisim'$ as above$,$   if
the associated  symplectic spreads  are equivalent  then
$n'=n,$
$F_i'=F_i$ and
$\zeta '_i=\zeta _i ^\sigma$ for some
 $\sigma\in\Aut F\2$~and all~$i$.}~We 
use these  theorems   to prove the following result, which  implies 
 \thmref{Main Theorem}:%

\medskip  
\begin{theorem} 
\label{orthogonal spread}
For odd $m>1,$ even $q^m>8$ and $\ffiti$ as above$,$
\begin{enumerate} 
\item[\rm(i)]  $\displaystyle
 \Sigma\ffiti:=\big\{\big(\sum _0^{n-1}W_i\vartheta_i \big)
\theta \mid \theta\in\c\big\}$ is an $O^-(2m,q)$-spread of the
 $F_{n}$-space $V$ equipped with the quadratic form $Q_n ,$
 and $\tc$ is a group of isometries acting  transitively on this elliptic spread$,$
and 
\item[\rm(ii)] $\Sigma\ffiti$  and  $\Sigma\ffisi'$ are equivalent if and only if 
$n'=n,$
$F_i'=F_i$ and
$\zeta '_i=\zeta _i ^\sigma$ for some $\sigma\in\Aut F\2$ and all $i$.

\end{enumerate} \rm
\end{theorem} 
  
   \proof
 (i) By \exref{symplectic to orthogonal} together with the above theorems, it suffices to verify that 
 the hyperplane 
 $\sum _0^{n -1}W_i\vartheta_i  $ 
 of  the t.i.~$F_n$-space 
$ \sum _0^{n -1}W_i\vartheta_i+F_{n }\vartheta_{n }  $ 
is its set of
 $Q_{n}$-singular vectors.
If   $w_i\in W_i ,  $  $0\le i\le m-1,$ then  
 $ Q_{n}(w_i\gamma_i  ) =   T_{n}(w_i\gamma_i 
 \overline{w_i\gamma_i })=
  T_{n} (w_i^2)=
 T_{n}T_{i+1}(w_i)^2=T_{n} (0)$.   Thus, 
  $\sum _0^{n -1}W_i\vartheta_i  $ is t.s. since it is a t.i.~subspace spanned by singular vectors.
   
  \smallskip
 (ii)  Assume that 
 $\Sigma\ffiti$  and  $\Sigma\ffisi'$ are equivalent by a semilinear map of $F\2$ preserving the orthogonal geometry.  Then $F\2$
 produces isomorphic vector spaces over its subfields
  $F_{n }$ and $F_{n '}$,  so that  $F_{n }=F_{n '}$
  and we have only one 
    quadratic form  
 $Q_{n } $ to consider.~(This avoids the    additional trace map
 $F\2\to \F_2$  used 
 in \cite[p.~8]{KW1}.)

The   cyclic group $\tc$   is transitive on both $\Sigma\ffisi' $  and $ \Sigma\ffiti$, so that   $\Sigma\ffisi' \!= \!\Sigma\ffiti^\sigma$ for some 
   $\sigma\in \Aut F\2$   by the exact same Sylow argument  as~in   \cite[proof of Theorem~5.2]{KW1}.  
   This equality of sets implies that 
   $n'=n$ and $\z_i'=\z_i^\sigma$     for some $\sigma\in\Aut F\2$ and
   all $i$
   as
 in  \cite[Lemma~5.3]{KW1} (more precisely, the case $l=l'=0$ of the proof
 of that lemma  is exactly  what is needed here).

The converse is trivial.
\qed
 
 \Remark\rm
As in \cite{KW1}, the classification of the finite simple groups was   not 
 needed for dealing with     elliptic spread equivalence.
However, in view of  \cite{BP}   the automorphism group of 
 $\Sigma\ffiti$ is a 1-dimensional semilinear group:   
 the semidirect product of  $\tc \times F_n^*$ with  the stabilizer  in 
 $\Aut F\2$
 of 
 all $\z_i$.     

\Remarks\rm
A symplectic spread in an $O^-(2n,q)$-space produces 
an elliptic spread (\exref{symplectic to orthogonal}).  
However, comparing  \eqnref{symplectic spread} and the results following it with
\thmref{orthogonal spread} shows that the same elliptic spread 
can arise from  many inequivalent symplectic spreads
by using different choices for $\z_n$.

This observation should be compared with \cite{Ka} and various sequels (such as 
\cite{KW1,KW2}).  Those papers  are based on the fact that 
a symplectic spread in 
an $\Sp(2m,q)$-space (with $m$ odd and $q$ even)  produces an essentially unique orthogonal  spread 
in an $O^+(2m+2,q)$-space, while an orthogonal  spread 
in an $O^+(2m+2,q)$-space produces many  inequivalent 
symplectic spreads and hence many affine  planes.%

\Remark\rm
Is there any way to decide whether or not a given 
elliptic spread $\Sigma$  in an $O^-(2n,q)$-space ``extends'' to   a symplectic spread?  In other words, is there a way to know from 
  properties of $\Sigma$ that, for each 
$X\in \Sigma$, it is possible to choose a t.i. subspace $\hat X\supset X$ so  that 
the set of these $\hat X$ is a symplectic  spread?

We used the cyclic group $\tc$ as a crutch for this purpose.
Such choices  presumably do  not arise   for the large numbers of  
$O^-(6,q)$-spreads obtained from a ``derivation'' process    \cite{CK,CE}. 
Moreover, although we deal with  large dimensions, the number of elliptic spreads we obtain in an  orthogonal  space  of 
fixed dimension is tiny compared to the number obtained by derivation in 
an $O^-(6,q)$-space.

\smallskip\smallskip
{\noindent  \bf Corrections.  1.}  Alan Prince has observed that   \cite[Theorem~1.1]{KW1} needs to be modified slightly by deleting the word  ``nondesarguesian''.   As it stands, that theorem states that there are more than $5/4$ nondesarguesian flag-transitive planes of
order 64, whereas  there is only one such plane \cite{Pr}
(constructed in \cite{KW1}).
\smallskip
\smallskip  

{\noindent  \bf 2.} In   \cite[Remark 6.6]{KW1} it   states that the affine planes obtained in that paper are \emph{precisely the flag-transitive
scions of the desarguesian plane of order $q^m$};
here ``scions'' refer to planes obtained by a recursive   ``up and down process'' described in \cite{KW1,KW2}.  That Remark  should also have continued
with the    assumption  that  these scions were   \emph{obtained by 
retaining flag-transitivity
throughout  the up and down process.}    Retaining flag-transitivity is needed for the inductive argument   in that Remark.  
Otherwise,   however unlikely it may seem,   this up and down process could magically produce a   spread having 
unexpected automorphisms. 

\smallskip   \smallskip
{\noindent\bf Acknowledgement.} \rm
This research was supported in part by a grant  
from the
Simons Foundation.


\begin{thebibliography}{GG1}
 
    \bibitem[BP]{BP}
   J. Bamberg and T. Penttila,
   A classification of transitive ovoids, spreads, and  $m$-systems
of polar spaces.
   Forum Math. 21 (2009)  181Ð-216.
   
     \bibitem[CE]{CE}
  A.  Cossidente and  G. L. Ebert,
Permutable polarities and a class of ovoids of the Hermitian surface. Eur. J. Comb. 25 (2004)  1059--1066. 


   
   \bibitem[CK]{CK}
   A. 
Cossidente and G. Korchm\'aros,  
Transitive ovoids of the Hermitian surface of $ \PG(3,q^2)$, $q$ even.  
J. Comb. Theory(A)  101 (2003)  117--130.    
 
  \bibitem[Di]{Dillon}
  J.~F. Dillon,
 Elementary Hadamard difference sets.
  Ph.D. thesis, U. of Maryland 1974.

   \bibitem[Dye]{Dye}   
R. H.   Dye, Partitions and their stabilizers for line complexes and quadrics.  Ann. Mat.
Pura Appl.  114 (1977)  173--194.

 \bibitem[JJB]{JJB}  
N. L. Johnson, V. Jha,   and M.  Biliotti,  Handbook of finite translation planes.  Chapman \& Hall/CRC, Boca Raton, FL  2007. 

   \bibitem[Ka]{Ka} W.~M.~Kantor,
Spreads, translation
planes and Kerdock sets. I, II. SIAM J. Algebraic and Discrete Methods 3   
(1982) 151--165,  308--318.

\bibitem[KW1]{KW1}$\!$W.~M.~Kantor and  M.~E.~Williams,
New flag-transitive affine  planes of even
order. J. Comb. Theory(A) 74 (1996) 1--13.

\bibitem[KW2]{KW2}$\!$W.~M.~Kantor and  M.~E.~Williams,
Symplectic semifield planes and $\Zsmall_4$-linear codes.  Trans. AMS  356 (2004)
895--938.

\bibitem[L\"u]{Lu}  H.
 L\"uneburg,   Die Suzukigruppen und ihre Geometrien.  Springer, Berlin-New York 1965
 

\bibitem[Pr]{Pr}
A. R. Prince,  Flag-transitive affine planes of order 
64.~Des. Codes Crypt.~18~(1999)~{{217--221}}.%

\bibitem[Ta]{Taylor} D. E.  Taylor, The geometry of the classical groups.
Heldermann, Berlin   1992. 

 
 
\end{thebibliography}
\end{document}